\documentclass[oneside,a4paper, 12pt]{article}
\usepackage{anysize}
\usepackage{long-geodesic}

\hypersetup{pdftitle={Long geodesics on convex surfaces},
pdfauthor={Arseniy Akopyan and Anton Petrunin}}

\begin{document}

\title{Long geodesics on convex surfaces}
\author{Arseniy Akopyan
\and 
Anton Petrunin}
\newcommand{\Addresses}{{\bigskip\footnotesize

Arseniy Akopyan, 
\par\nopagebreak\textsc{Institute of Science and Technology Austria,} 
\par\nopagebreak\textsc{Am Campus 1, 3400 Klosterneuburg, Austria.}
\par\nopagebreak\textit{Email}: \texttt{akopjan@gmail.com}

\medskip

Anton Petrunin, 
\par\nopagebreak\textsc{Mathematics Department,} 
\par\nopagebreak\textsc{Pennsylvania State University,}
\par\nopagebreak\textsc{University Park, PA 16802, USA}
  \par\nopagebreak
  \textit{Email}: \texttt{petrunin@math.psu.edu}
}}
\date{}
\maketitle

\begin{abstract}
We review the theory of intrinsic geometry of convex surfaces in the Euclidean space and prove the following theorem: 
if the boundary surface of a convex body $K$ contains arbitrarily long closed simple geo\-de\-sics, then $K$ is an isosceles tetrahedron.
\end{abstract}


\section{Introduction}
The goal of this note is to introduce the reader to the theory of intrinsic geometry of convex surfaces.
We illustrate the power of the tools by proving a theorem on convex surfaces containing an arbitrarily long closed simple geodesic. 

Let us remind that a curve in a surface is called \emph{geodesic} if every sufficiently short arc of the curve is length minimizing; if in addition it has no self intersections, we call it \emph{simple geodesic}.
A tetrahedron with equal opposite edges will be called \emph{isosceles}.

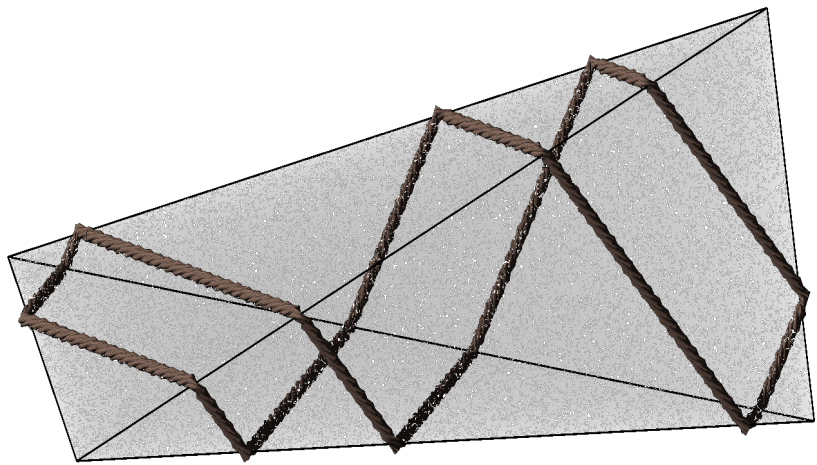
\begin{center}
\begin{lpic}[t(-1 mm),b(-1 mm),r(0 mm),l(0 mm)]{pics/snow(1)}
\end{lpic}
\end{center}

\begin{thm}{Theorem}\label{Long geodesic}
Assume that the boundary surface $\Sigma$ of a convex body $K$ in the Euclidean space $\EE^3$
admits an arbitrarily long simple closed geodesic.
Then $K$ is an isosceles tetrahedron.
\end{thm}

This gives an affirmative answer to the question asked by Vladimir Protasov;
he proved the statement in assumption that $\Sigma$ is the boundary surface of a convex polyhedron, see \cite{protasov2008onthenumber, protasov2012}.

The axiomatic method of Alexandrov geometry allows to work with metrics of convex surfaces directly, without approximating it first by smooth or polyhedral metric.
Such approximations destroy the closed geodesics on the surface; therefore
it is hard (if at all possible) to apply approximations in the proof of our theorem.
On the other hand, a proof in the smooth or polyhedral case usually admits a \emph{translation} into \emph{Alexandrov's language};
such translation makes the result more general.
In fact, our proof resembles a translation of the proof given by Protasov.

Note that the main theorem implies in particular that a smooth convex surface does not have arbitrarily long simple closed geodesics.
However we do not know a proof of this corollary which is essentially simpler than the one presented below.

\section{An overview of the theory}

\begin{wrapfigure}[4]{r}{21 mm}
\begin{lpic}[t(-6 mm),b(0 mm),r(0 mm),l(0 mm)]{pics/cosine-rule(1)}
\lbl[t]{10,1;$a$}
\lbl[br]{5.5,9.5;$b$}
\lbl[bl]{16,9;$c$}
\lbl[bl]{7,4;$\gamma$}
\end{lpic}
\end{wrapfigure}

The geometry of intrinsic metric on convex surfaces
can be considered as generalization of Euclidean plane geometry,
where instead of equality in the cosine rule we have the following inequality
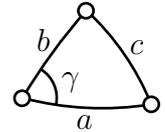
\[c^2\le a^2+b^2-2\cdot a \cdot b\cdot \cos\gamma.\]
This approach provides a uniform way to work with smooth and singular surfaces.

As the main reference, we will use the classical  book of Alexandrov \cite{aleksandrov1948vnutrennnyaya}; it is available in Russian and German.
The books \cite{busemann1958convex} or \cite{pogorelov1973extrinsic} of Busemann and Pogorelov should work as well and they are available in English.

As a motivation, we recommend an other classical book of Alexandrov \cite{aleksandrov2005} --- it is fun to read, it worth reading and once red, it should be easy to read it should also be easy to read most papers on the subject including \cite{aleksandrov1948vnutrennnyaya}.
Some related topics discussed in popular books; see for example \cite[Lectures 20, 24 and 25]{tabachnikov-fuks}.

\subsection*{Angles}

By a \emph{surface} we mean a compact $2$-dimensional manifold
(possibly with non-empty boundary), which is equipped with \emph{geodesic} metric.
A metric on the surface $\Sigma$ is called geodesic if any two points $p,q\in \Sigma$ can be joined by a curve with length $|p-q|_\Sigma$;
we denote by $|p-q|_\Sigma$ the distance from $p$ to $q$ in $\Sigma$.
This curve will be denoted by $[pq]$ and called a \emph{minimizing geodesic} from $p$ to~$q$.

A triple of points $x,y,z\in\Sigma$, together with a choice of three minimizing geodesics $[x y]$, $[y z]$, $[z x]$ will be called a \index{triangle}\emph{triangle} 
and denoted by $[x y z]$.
Many different triangles with vertices $x$, $y$ and $z$ may exist, 
any of which can be denoted by $[x y z]$.

\begin{wrapfigure}{r}{19 mm}
\begin{lpic}[t(-0 mm),b(0 mm),r(0 mm),l(0 mm)]{pics/hinge(1)}
\lbl[t]{17,2;$z$}
\lbl[t]{8,0;$\bar z$}
\lbl[r]{5,17;$y$}
\lbl[r]{1.5,8;$\bar y$}
\lbl[rt]{1,0;$x$}
\end{lpic}
\end{wrapfigure}

A triangle $[\~x\~y\~z]$ in the Euclidean plane $\EE^2$
with the same side lengths as $[x y z]$ 
is called a \emph{model triangle} of $[x y z]$;
this relation will be written as $[\~x\~y\~z]=\~\triangle(x y z)$.
The angle $\mangle\hinge{\~x}{\~y}{\~z}$ of the model triangle $[\~x\~y\~z]$ is called \emph{model angle} of the triangle $[x y z]$ at $x$ and denoted by $\angk x y z$.

A pair of geodesics $[xy]$ and $[xz]$ with common endpoint $x$ is called a \emph{hinge} and is denoted by $\hinge{x}{y}{z}$.
The angle measure $\mangle\hinge{x}{y}{z}$ of the hinge is defined as the limit of model angles for the triangles sliding along the sides of hinge to its vertex. 
That is,
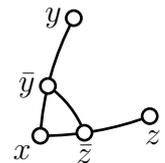
\[\mangle\hinge{x}{y}{z}=\lim_{\bar y,\bar z\to x}\set{\angk{x}{\bar y}{\bar z}}{ \bar y\in \left]xy\right], \bar z\in \left]xz\right]},\]
where $\left]xy\right]=[xy]\backslash\{x\}$.
In general, the angle of a hinge maybe undefined, but as you will see it will be defined every time we need it.
Note that if $p\in \left]xy\right[=[xy]\backslash\{x,y\}$ then $\mangle\hinge{p}{x}{y}=\pi$.

\subsection*{Comparison}

\begin{thm}{Comparison property}\label{Comparison property}
We say that a hinge $\hinge x y z$ 
satisfies the \emph{comparison property} if the angle
$\mangle\hinge{x}{y}{z}$ is defined and 
\[\mangle \hinge{x}{y}{z} \ge \angk{x}{y}{z}{}{}.\]
\end{thm}

Two hinges $\hinge{p}{x}{z}$ and $\hinge{p}{y}{z}$ will be called \emph{supplementary},
if they share a side $[pz]$ and $p\in \left]xy\right[$.

\begin{wrapfigure}{r}{27 mm}
\begin{lpic}[t(-0 mm),b(0 mm),r(0 mm),l(0 mm)]{pics/suplimentary(1)}
\lbl[t]{24,2.6;$x$}
\lbl[t]{1,13.5;$y$}
\lbl[t]{5,-.5;$z$}
\lbl[bl]{15,13;$p$}
\end{lpic}
\end{wrapfigure}

\begin{thm}{Supplementary property}\label{Supplementary property}
We say that the \emph{suplementary property} holds for two supplementary hinges $\hinge{p}{x}{z}$ and $\hinge{p}{y}{z}$ if the angles $\mangle\hinge{p}{x}{z}$ and $\mangle\hinge{p}{y}{z}$ are defined and
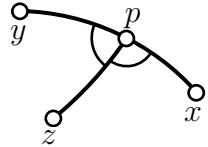
\[\mangle\hinge{p}{x}{z}+\mangle\hinge{p}{y}{z}=\pi.\]

\end{thm}

We will say that a surface $\Sigma$ has \emph{non-negative curvature in the sense of Alexandrov}
if the angles of all hinges in $\Sigma$ are defined and satisfy the comparison and supplementary properties.%
\footnote{It is not known, if the comparison property for all hinges in $\Sigma$ imply the supplementary property.}

The geometry of such surfaces is very specific. 
For example, by the comparison property, if a hinge has vanishing angle then one of its sides lies in the other.
In particular, geodesics in $\Sigma$ can not bifurcate. The following theorem plays the central role in the theory.

\begin{thm}{Globalization theorem}\label{Globalization theorem}
Assume that any point of the surface $\Sigma$ has a neighborhood $U$ such that the angles of all hinges in $U$ are defined and satisfy 
the comparison and supplementary properties.
Then these properties hold for all hinges;
that is, $\Sigma$ has non-negative curvature in the sense of Alexandrov.
\end{thm}

This theorem was proved by Alexandrov \cite{alexandrow1957ubereine} 
and generalized since then many times;
an amusing proof was given recently by Urs Lang and Viktor Schroeder in \cite{lang2012toponogov}.

Further, assume that the surface $\Sigma$ has non-negative curvature in the sense of Alexandrov,
$\gamma$ is a geodesic in $\Sigma$ parametrized by length, 
and let $p$ be an arbitrary point on $\Sigma$.
Note that the function $f(t)=|p-\gamma(t)|_\Sigma$ is  $1$-Lipschitz. 
In particular, the function $f$ is differentiable almost everywhere.

\begin{wrapfigure}{r}{25 mm}
\begin{lpic}[t(-0 mm),b(0 mm),r(0 mm),l(0 mm)]{pics/first-variation(1)}
\lbl[t]{5,0;$p$}
\lbl[bl]{15,13;$\gamma(t)$}
\lbl[rt]{10,12;$\phi_{+}$}
\lbl[lt]{15,8;$\phi_{-}$}
\end{lpic}
\end{wrapfigure}

Let us denote by $\phi_\pm(t)$ the angles between the positive and negative directions of $\gamma$ and a geodesic $[\gamma(t)\,p]$; 
see the diagram.
From the definition of angle, via the triangle inequality,
one gets the following
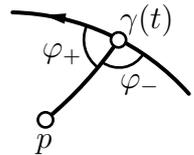
\[\pm f'(t)\le -\cos[\phi_\pm(t)]\]
at any $t$ for which the derivative $f'(t)$ is defined (see  \cite[XI \S 2 (7)]{ aleksandrov1948vnutrennnyaya}).

By the supplementary property $\phi_-\z+\phi_+=\pi$;
hence $\cos\phi_+ +\cos\phi_-=0$.
Therefore the two inequalities above imply so called \emph{first variation formula}
\begin{equation}
	\label{eq:first variation}
f'(t)=-\cos[\phi_+(t)]
	\tag{${*}$}
\end{equation}
for any $t$, where the derivative $f'(t)$ is defined.

Further, for any surface $\Sigma$ with non-negative curvature in the sense of Alexandrov,
the Kirszbraun extension theorem holds.
That is, any distance non-expanding map from a subset of $\Sigma$ to the Euclidean plane can be extended to a distance non-expanding map defined on whole $\Sigma$;
see \cite{lang1997kirszbraun, alexander2011alexandrov}.
(In fact, this statement could be used to define the spaces with non-negative curvature in the sense of Alexandrov.)

Applying the Kiszbraun theorem for three-point sets,
one gets the following area comparison property which will be important to us.

\begin{thm}{Area comparison}\label{Area comparison}
Assume $\Sigma$ is a surface with non-negative curvature in the sense of Alexandrov
and a triangle $[xyz]$ in $\Sigma$ bounds a open set $\Delta$ homeomorphic to a disc.
Then 
\[\area\Delta\ge \area\~\triangle(xyz).\]
\end{thm}
An earlier proof of this theorem by slicing $\Delta$ into small triangles is given in \cite[X \S 1]{ aleksandrov1948vnutrennnyaya}.

\subsection*{Convex surfaces}

Recall that the intrinsic distance between points $x$ and $y$ on a surface $\Sigma$ in $\EE^3$, is defined as the  greatest lower bound for the lengths of curves connecting $x$ to $y$ in $\Sigma$.

\begin{thm}{Comparison theorem}\label{Comparison theorem}
The boundary surface of any convex body in $\EE^3$,
if equipped with the induced intrinsic metric, 
is a sphere with a non-negatively curved metric in the sense of Alexandrov.
\end{thm}

The converse of this theorem also holds
if one considers convex plane figure as a degenerate convex body \cite[III \S 3]{aleksandrov1948vnutrennnyaya}.
The \emph{boundary surface} of a flat convex figure has to be defined as its doubling;
that is, two copies 
of the figure glued along the boundary --- 
it will look like a surface if you can walk on both sides of the figure, but can not pass from one side directly to an other side.

\begin{thm}{Theorem}
Any surface $\Sigma$ with non-negative curvature in the sense of Alexandrov which is homeomorphic to the sphere
is isometric to the boundary surface of a convex body in $\EE^3$,
which possibly degenerates to a flat figure.
\end{thm}

It turns out that $\Sigma$ defines the convex body up to congruence.
This is a very hard theorem;
it was first proved by Pogorelov in \cite{pogorelov1952odnoznacnaya}.
We will use the following weaker statement proved by Alexandrov earlier 
\cite[VI \S 5]{aleksandrov1948vnutrennnyaya}, \cite[3.3]{aleksandrov2005};
its proof goes essentially the same way as Cauchy's rigidity theorem.

\begin{thm}{Rigidity theorem}\label{Rigidity theorem}
Convex polyhedrons in $\EE^3$ with isometric boundary surfaces are congruent. 
\end{thm}

\subsection*{Curvature measure}

Let $\Sigma$ be a surface with non-negative curvature in the sense of Alexandrov.
Assume a triangle $[xyz]$ bounds an open set $\Delta$ in $\Sigma$ which is convex and homeomorphic to a disc.
By convex we mean that any minimizing geodesic with endpoints in $\Delta$ lie completely in $\Delta$.
In this case we define the \emph{curvature} of $\Delta$ as the \emph{angle excess} of the triangle $[xyz]$;
that is,
\[\kappa(\Delta)=\mangle \hinge x y z+\mangle \hinge  y z x+\mangle \hinge z x y-\pi.\]
The functional $\kappa$ uniquely extends to non-negative measure, the so called \index{curvature measure}\emph{curvature measure}, defined on all Borel subsets in the interior of $\Sigma$.

\begin{wrapfigure}{r}{23 mm}
\begin{lpic}[t(-4 mm),b(-0 mm),r(0 mm),l(0 mm)]{pics/geodesic(1)}
\end{lpic}
\end{wrapfigure}

It turns out that any geodesic without its end-points has vanishing curvature.
This can be proved by covering the interior of geodesic by two thin triangles with small excess as shown in the picture.
(This requires some work, but simple.)

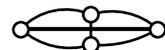
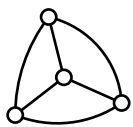
\begin{wrapfigure}{l}{20 mm}
\begin{lpic}[t(-0 mm),b(-0 mm),r(0 mm),l(0 mm)]{pics/mercedes(1)}
\end{lpic}
\end{wrapfigure}

Further, the curvature of an interior point in $\Sigma$ is $2{\cdot}\pi$ minus total angle around it.
The later can be seen from the picture ---
to find the curvature of the central vertex one has to subtract excesses of three small triangles from the excess of the big one. (The existence of such configuration also requires some work.)
A point of non-zero curvature is called \emph{singular}.

For curvature measure, an analog of the Gauss--Bonnet formula holds;
in particular, if $\Sigma$ is a surface with non-negative curvature in the sence of Alexandrov then
\begin{itemize}
\item If $\Sigma$ is homeomorphic to the sphere then 
\[\kappa(\Sigma)=4{\cdot}\pi\]
\item If a closed geodesic cuts from $\Sigma$ a disc $\Delta$ then 
\[\kappa(\Delta)=2{\cdot}\pi.\]
\item If a closed broken geodesic cuts from $\Sigma$ a disc $\Delta$ then 
\[\kappa(\Delta)=2{\cdot}\pi +(\theta_1-\pi)+\dots+(\theta_n-\pi),\]
where $\theta_1,\dots,\theta_n$ are the inner angles at the corners of $\Delta$.
\end{itemize}

\section{Four singular points}

The following lemma is the key to the proof of Theorem~\ref{Long geodesic}.
A point in a surface is called \emph{flat} if it admits a flat neighborhood;
that is, a neighborhood isometric to an open subset of the plane.
Equivalently, the curvature measure is vanishing in a neighborhood of this point.
Note that for a polyhedral surface all points except the vertexes are flat.

\begin{thm}{Lemma} 
	\label{lem:4 singular points}
Assume that the boundary surface $\Sigma$ of a convex body $K$ in $\EE^3$
admits an arbitrarily long simple closed geodesic.
Then $\Sigma$ contains $4$ singular points with curvature $\pi$ and the rest of it is flat.
\end{thm}

In the proof we will show that 
it is possible to find four non-intersecting open sets of arbitrarily small diameter, 
such that each of them has curvature arbitrary close to $\pi$.
Passing to the limit we get the statement of the lemma.

By cutting the surface $\Sigma$ along a sufficiently long closed simple geodesic,
we get two discs.
The key step is to show that each of these discs is ``long'' and ``thin''.
Then we show that their ``ends'' form four required sets with big curvature.
To be more precise, the four ends will have small perimeter,
but according to the following claim, it implies small diameter as well.
In fact, this claim holds for any metric on the sphere, not necessary the closed convex surface.

\begin{thm}{Claim}\label{Lemma:diameter-perimeter}
For any $\varepsilon>0$ there is a $\delta>0$, such that any simple curve in $\Sigma$ of length smaller than $\delta$ bounds a region of diameter at most $\varepsilon$.
\end{thm}

\parit{Proof.}
We will come to contradiction by showing that there is a point which cuts $\Sigma$ into connected non-empty parts.
Suppose there is a sequence of curves $\gamma_n$ which cuts $\Sigma$ into open regions $A_n$ and $B_n$ of diameter at least $\varepsilon$ each and such that $\length\gamma_n\to 0$ as $n\to\infty$. 
By compactness of $\Sigma$,
we can pass to a subsequence of $\gamma_n$ which converges to a point, denote it by $p$.

Note that for large $n$, each of the regions $A_n$ and $B_n$ contains a disk of radius $\tfrac{\varepsilon}{3}$;
label their centers by $a_n$ and $b_n$ respectively. 
Indeed, if $n$ is large, then $\gamma_n$ lies in $\tfrac{\varepsilon}{6}$ neighborhood of $p$.
Take $a_n$ and $b_n$ to be points in the regions maximizing the distance to $p$, note that these distances should be greater than $\tfrac{\varepsilon}{2}$.
Now, note that if disk of radius $\tfrac{\varepsilon}{3}$ centered at $a_n$ intersects $\gamma_n$, then distance from $a_n$ to $p$ will be less then $\tfrac{\varepsilon}{3}+\tfrac{\varepsilon}{6}=\tfrac{\varepsilon}{2}$.

Pass to a subsequence of $\gamma_n$ so that $a_n$ and $b_n$ converge, denote by $a$ and $b$ their limits.
Note that for large $n$ the domains $A_n$ and $B_n$ contain the disks of radius $\tfrac\varepsilon4$ centered at $a$ and $b$ correspondingly.

Since for all large $n$, any path from $a$ to $b$ has 
to cross $\gamma_n$ and the sequence $\gamma_n$ converges to~$p$, we obtain that any path from $a$ to $b$ contains $p$.
That is, $p$ is a cut point of~$\Sigma$, a contradiction.
\qeds

\parit{Proof of Lemma~\ref{lem:4 singular points}.}
According to the comparison theorem (\ref{Comparison theorem}) the surface $\Sigma$ has non-negative curvature in the sense of Alexandrov.
Cut $\Sigma$ along a sufficiently long closed simple geodesic,
we get two discs.
Choose one of the discs, say $D$. 
Note that $D$ is locally convex in $\Sigma$; that is, sufficiently short minimizing geodesic with the ends in $D$ has to lie in $D$.

Equip $D$ with the intrinsic metric further denoted by $|{*}-{*}|_D$.
By the globalization theorem (\ref{Globalization theorem}), $D$ has non-negative curvature in the sense of Alexandrov.

{

\begin{wrapfigure}{r}{31mm}
\begin{lpic}[t(0 mm),b(-0 mm),r(0 mm),l(1 mm)]{pics/long-geodesic-diam(1)}
\lbl[r]{0,8;$p$}
\lbl[l]{28,8;$q$}
\lbl[b]{10,16;$\gamma_1(t)$}
\lbl[l]{12,9,-60;$\ge\tfrac\pi3$}
\end{lpic}
\end{wrapfigure}

Choose a pair of points $p,q\in\partial D$ which maximize the distance $|p-q|_D$.
Clearly,
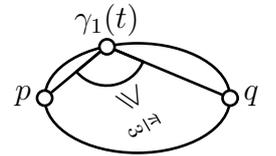
\[|x-p|_D\le |p-q|_D,\quad |x-q|_D\le |p-q|_D\] 
for any other point $x\in\partial D$.
By the comparison property (\ref{Comparison property}) 
\begin{equation}
	\label{eq:pxq>pi/3}
	\measuredangle[x\,^p_q]\ge \tfrac\pi3.
	\tag{${**}$}
\end{equation}

The points $p$ and $q$ divide $\partial D$ into two arcs,
say $\gamma_1$ and $\gamma_2$;
let us parametrize them by arc length from $p$ to $q$. 
By the first variation formula \eqref{eq:first variation}, for almost all $t$ we have
\[\tfrac{d}{dt}|p-\gamma_i(t)|_D=-\cos \phi,\] 
where $\phi$ is the angle between the direction of $\gamma_i$ at $x=\gamma_i(t)$ and the geodesic $[xp]$.
The same way 
\[\tfrac{d}{dt}|q-\gamma_i(t)|_D=-\cos \psi,\] 
where $\psi$ is the angle between the direction of $\gamma_i$ at $x\z=\gamma_i(t)$ and the geodesic $[xq]$.
By \eqref{eq:pxq>pi/3} and the supplementary property (\ref{Supplementary property}) we have $\phi-\psi\ge \tfrac\pi3$.
Therefore 
\begin{equation*}
\tfrac{d}{dt}\left(|p-\gamma_i(t)|_D-|q-\gamma_i(t)|_D\right)
= \cos \psi-\cos\phi
\ge
\tfrac12.
\end{equation*}
The function $t\mapsto|p-\gamma_i(t)|_D-|q-\gamma_i(t)|_D$ ranges from $-|p-q|_D$ to $|p-q|_D$.
Therefore
\[|p-q|_D\ge \tfrac14{\cdot}\length\gamma_i\]
and therefore
\[|p-q|_D\ge \tfrac18{\cdot}\length[\partial D].\tag{$\asterism$}\]
That is, if the geodesic was long 
then $D$ has large diameter.

}

Fix small $\eps>0$ and a positive $\delta\ll \eps$. 
(The value $\delta=\varepsilon\cdot\sin \varepsilon$ will do.
It means that in any Euclidean triangle with one side $\varepsilon$ and another $\delta$ the angle opposite to the side of length $\delta$ is less than~$\varepsilon$.)

Let $x$ and $y$ be the first points along $\gamma_1$ and $\gamma_2$ respectively at distance $\varepsilon$ from $p$.
If $|q-x|_D>|q-y|_D$, move the point $y$ along $\gamma_2$ toward $p$ until $|q-x|_D=|q-y|_D$; 
since $|p-q|_D$ is maximal, this distance will be achieved. 
In case $|q-x|_D<|q-y|_D$ do the same for the point $x$. 
Now we have
$$|q-x|_D=|q-y|_D\ge|p-q|_D-\varepsilon.$$

By the area comparison (\ref{Area comparison})
\[\area [\tilde \triangle qxy] \le \area [\triangle qxy] \le \area \Sigma.\]
Since $\tilde \triangle qxy$ is isosceles, its area is larger than $\tfrac1{10}\cdot|x-y|_D\cdot|x-q|_D$;
this holds since $\eps$ is small.
Applying $(\asterism)$, we get
\begin{equation*}
\label{eq:|x-y|}
\begin{aligned}
|x-y|_D&\le
100\cdot\frac{ \area\Sigma}{\length[\partial D]}.
\end{aligned}
\end{equation*}
Therefore, if $D$ has large perimeter then we can assume that $|x-y|_D <\delta$.

Cut $D$ by $[xy]$
and consider the part (a lune) $L_p$ with the point $p$ in it.
Note that the curvature of $L_p$ is $\alpha+\beta$, where $\alpha$ and $\beta$ are the angles as on the diagram.
By the comparison property $\alpha\ge \tilde\measuredangle(x\,^p_y)$ 
and $\beta\ge \tilde\measuredangle(y\,^p_x)$.
By Gauss--Bonnet formula, the curvature of $L_p$ is at least $\pi-\tilde\measuredangle(p\,^x_y)>\pi-\varepsilon$.

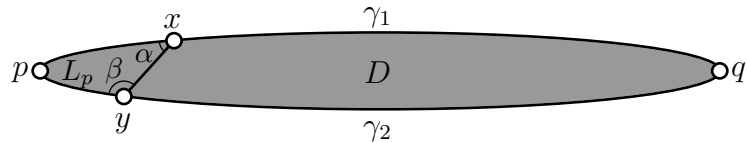
\begin{center}
\begin{lpic}[t(3 mm),b(3 mm),r(0 mm),l(0 mm)]{pics/long-geodesic-D(1)}
\lbl{46,6;$D$}
\lbl{6.5,5.7;{\small $L_p$}}
\lbl[r]{0,6;$p$}
\lbl[l]{92.5,6;$q$}
\lbl[b]{19,11.8;$x$}
\lbl[t]{12.5,.5;$y$}
\lbl[b]{46,12;$\gamma_1$}
\lbl[t]{46,-.5;$\gamma_2$}
\lbl[tr]{16.5,9;{\small $\alpha$}}
\lbl[br]{12.5,4.5;{\small $\beta$}}
\end{lpic}
\end{center}

Assuming $\partial D$ is long, we can find a lune $L_p$ with perimeter at most $2\cdot\eps + \delta$,
such that curvature $L_p$ is at least $\pi-\eps$.
If $\eps$ is small, by Lemma \ref{Lemma:diameter-perimeter}, $L_p$ has small diameter, say $\eps'=\eps'(\eps)$.

Using the same construction for $p$ and $q$ in the disc $D$,
and for the other disc,
we get four lunes in $\Sigma$, 
each of diameter at most $\eps'$, 
and each with curvature at least $\pi-\eps$.
By Gauss--Bonnet formula, the remaining curvature is at most $4\cdot\eps$.

Since $\eps>0$ is arbitrary, it follows that support of curvature measure can be covered by four sets with arbitrarily small diameter;
that is, support of $\kappa$ is a $4$-point set.
Clearly each of these points has curvature $\pi$,
the remaining part of $\Sigma$ has vanishing curvature and therefore flat.
\qeds

\section{Isosceles tetrahedron}

\begin{thm}{Lemma} 
Assume that a closed convex surface in $\EE^3$
has $4$ singular points with curvature $\pi$.
Then it bounds an isosceles tetrahedron.
\end{thm}

{

\begin{wrapfigure}{o}{21 mm}
\begin{lpic}[t(-4 mm),b(-3 mm),r(0 mm),l(0 mm)]{pics/akopyan(1)}
\lbl[br]{3,8;$x$}
\lbl[bl]{14,8;$z$}
\lbl[t]{10,0;$y$}
\end{lpic}
\end{wrapfigure}

\parit{Proof.}
Denote the surface by $\Sigma$ and by $p$, $x$, $y$ and $z$ the singular points.

Connect $p$ to the $x$, $y$ and $z$ by three minimizing geodesics $[px]$, $[py]$ and $[pz]$.
Note that these geodesics intersect only at the common point $p$.

Let us cut $\Sigma$ along three geodesics $[px]$, $[py]$ and $[pz]$.
Develop the obtained flat surface on the plane.
Since each point $x$, $y$ and $z$ have curvature $\pi$,
the pairs of cuts at these points open straight.
Therefore, the development forms a triangle; denote it by $\triangle$.
The points $x$, $y$ and $z$ correspond to the midpoints of the sides of $\triangle$
and the point $p$ correspond to the vertices of $\triangle$.

}

It follows that $\Sigma$ is isometric to the boundary surface of isosceles tetrahedron.
The statement now follows from the Pogorelov's uniqueness theorem,
mentioned in the theory overview.

Alternatively one could apply the following exercise and the rigidity theorem (\ref{Rigidity theorem}) which has a much simpler proof.
\qeds

\begin{thm}{Exercise}
Assume that a convex body $K$ is bounded by a closed flat surface with finite number of singular points.
Show that $K$ is a polyhedron.
\end{thm}

\bigskip
\textbf{Acknowledgment.}
We want to thank Kristof Huszar and Sergei Tabachnikov for valuable comments on the preliminary version of this note.

\Addresses

\end{document}